\documentclass[11pt,letterpaper,reqno]{amsart}

\usepackage{amsfonts}
\usepackage{amsmath}
\usepackage{amsthm, amscd}
\usepackage{mathtools}
\usepackage{amssymb}
\usepackage{mathrsfs}
\usepackage{graphicx}
\usepackage{caption}
\usepackage{subcaption}
\usepackage{float}
\usepackage{xypic}
\usepackage[abs]{overpic}
\usepackage[alphabetic]{amsrefs}
\usepackage{etex}
\usepackage{tikz-cd}

\usepackage{hyperref}
\hypersetup{
	colorlinks,
	linkcolor=[rgb]{0,0,0.7},
	urlcolor=[rgb]{0,0,0.4},
	citecolor=[rgb]{0.4,0.1,0}
}

\allowdisplaybreaks

\theoremstyle{plain}\newtheorem{Theorem}{Theorem}[section]
\theoremstyle{plain}
\theoremstyle{plain}
\theoremstyle{plain}\newtheorem{Definition}[Theorem]{Definition}
\theoremstyle{plain}\newtheorem{Proposition}[Theorem]{Proposition}
\theoremstyle{plain}
\theoremstyle{plain}
\theoremstyle{plain}\newtheorem*{Claim*}{Claim}
\theoremstyle{plain}\newtheorem*{Theorem*}{Theorem}
\theoremstyle{plain}

\theoremstyle{remark}\newtheorem{remark}[Theorem]{Remark}
\theoremstyle{remark}\newtheorem{Example}[Theorem]{Example}
\theoremstyle{remark}\newtheorem*{Notation*}{Notation}

\theoremstyle{plain}
\makeatletter
\newtheorem*{rep@theorem}{\rep@title}
\newcommand{\newreptheorem}[2]{
\newenvironment{rep#1}[1]{
 \def\rep@title{#2 \ref{##1}}
 \begin{rep@theorem}}
 {\end{rep@theorem}}}
\makeatother
\newreptheorem{Theorem}{Theorem}
\newreptheorem{Proposition}{Proposition}
\newreptheorem{Corollary}{Corollary}

\numberwithin{equation}{section}

\usepackage{lipsum}

\DeclareMathOperator{\Diff}{Diff}
\DeclareMathOperator{\BDiff}{BDiff}
\DeclareMathOperator{\Homeo}{Homeo}
\DeclareMathOperator{\spinc}{\mathrm{spin}^\mathbb{C}}
\DeclareMathOperator{\Aut}{Aut}

\title{A note on the boundary Dehn twist of $K3$ surfaces}
\author{Yujie Lin}
\address{Qiuzhen College, Tsinghua University, Beijing 100084, China}
\email{lyj23@mails.tsinghua.edu.cn}

\begin{document}
\maketitle
\begin{abstract}
    By the work of Baraglia-Konno and Kronheimer-Mrowka, the boundary Dehn twist on punctured $K3$ surfaces is nontrivial in the smooth mapping class group relative to boundary. In this short note, we prove that it becomes trivial after abelianization. The proof is based on an obstruction for $\mathrm{Spin}^\mathbb{C}$ families due to Baraglia-Konno and the global Torelli theorem of $K3$ surfaces. 
\end{abstract}
\section{Introduction}
 Let $X$ be a smooth closed 4-manifold with $\pi_1(X)=1$. Denote $X^\circ=X-D^4$. Then $\partial X^\circ\cong S^3$ has a collar neighbourhood $N\cong S^3\times [0,1]$. Let $\gamma:[0,1]\to \mathrm{SO}(4)$ be a smooth non-contractible loop which is $id$ near $0$ and $1$ . Then the boundary Dehn twist $t_X$ of $X^\circ$ is defined as $t_X(x,t)=(\gamma_t(x),t)$ on $N\cong S^3\times[0,1]$ and extends by identity on the complement of $N$.  So $t_X\in \Diff_\partial(X^\circ)$ is a diffeomorphism of $X^\circ$ that fixes the boundary. 
 
 A natural and important question is whether $t_X$ is topologically or smoothly isotopic to identity relative to boundary. Here we only discuss the case $\pi_1(X)=1$. In the topological case, $[t_X]=1\in\Homeo_\partial(X^\circ)$ (see \cite{orson2023mapping}\cite{krannich2025torelli}). In the smooth case, Giansiracusa \cite{giansiracusa2008stable} proved that the kernel of the surjective homomorphism $\pi_0(\Diff_\partial(X^\circ))\to \pi_0(\Diff(X))$ is generated by $[t_X]$, hence the kernel is either 0 or $\mathbb{Z}/2$. It has been known that $[t_X]=1\in\Diff_\partial(X^\circ)$ for non-spin $X$ (see \cite{orson2023mapping}). For spin $X$, we don't have a general conclusion currently. But there are some partial results: 
\begin{itemize}
    \item (Giansiracusa \cite{giansiracusa2008stable}) $[t_X]=1$ for $X=\#_n(S^2\times S^2)$.
    \item (Baraglia-Konno \cite{baraglia2022bauer}, Kronheimer-Mrowka \cite{kronheimer2021dehn}) $[t_X]\neq 1$ for homotopy $K3$ surface.
    \item (Lin \cite{lin2023isotopy}) $[t_X]\neq 1$ for $X=K3\#(S^2\times S^2)$.
    \item (Baraglia-Konno \cite{baraglia2024irreducible}) $[t_X]\neq 1$ for those $X$ satisfying: $\pi_1(X)=1$, $b_+(X)\equiv 3 \mod 4$,  $\sigma(X)\equiv 16 \mod 32$, and there exists a $\spinc$ structure $\mathfrak{s}$ on $X$ such that $32\mid c_1(\mathfrak{s})$, $\mathrm{SW}(X,\mathfrak{s})$ is odd. For example, a large family of elliptic surfaces and complete intersections satisfy this technical condition.
\end{itemize}
 Cohomology classes of $\BDiff(X)$ correspond to characteristic classes of smooth $X-$bundle. So it's natural to study $H_*(\BDiff_\partial(X^\circ))$. Since $\pi_0(\Diff_\partial(X^\circ))\cong \pi_1(\BDiff_\partial(X^\circ))$, the abelianization of $\pi_0(\Diff_\partial(X^\circ))$ is canonically isomorphic to $H_1(\BDiff_\partial(X^\circ))$. A natural question is whether a nontrivial boundary Dehn twist $t_X$ becomes trivial under abelianization. In this short note we give an affirmative answer for $X=K3$.
 \begin{Theorem}\label{11}
For $X=K3$, $[t_X]_{ab}=0\in H_1(\BDiff_\partial(X^\circ))$
\end{Theorem}
Theorem \ref{11} follows from the following proposition.  
\begin{Proposition}\label{12}
    There exists a smooth fiber bundle $K3\to E\to T^2$ such that $E$ is not spin.
\end{Proposition}

There are two key ingredients in the proof of Proposition \ref{12}, both due to Baraglia-Konno \cite{baraglia2022bauer}\cite{baraglia2023note}. The first one (Theorem \ref{22}) is the relation between the first Chern class of the index family and the second Stiefel-Whitney class of the $H^+$ bundle for a $\spinc$ family of 4-manifold, which is proved by computing Steenrod squares of the mod 2 families Seiberg-Witten invariants. The second one (Theorem \ref{23}) is the existence of a section $s_{K3}$ of the natural homomorphism $\varphi_{K3}: \pi_0(\Diff(K3))\rightarrow \Aut(H^2(K3;\mathbb{Z}))$ over its image $\Gamma_{K3}$, which is proved using the global Torelli theorem for $K3$ surfaces.

\textbf{Outline of the proof.} In Section \ref{sec2}, we prove the following criterion: the triviality of the abelianized boundary Dehn twist $[t_X]_{ab}$ is equivalent to the existence of a smoooth $X$-bundle over a closed oriented surface without family spin structure. To construct such a bundle, in Section \ref{sec3}, we construct two specific commuting isomorphisms in $\Aut(K3;\mathbb{Z})$ and use the section $s_{K3}$ to realize them as two diffeomorphisms which commute up to isotopy. This gives a smooth $K3$-bundle on the torus with $w_2(H^+)\neq 0$. Suppose it admits a family spin structure, then the first Chern class of the index family is zero, which is impossible by the relation between index family and $H^+$ bundle.

\textbf{Acknowledgements.} The author is grateful to his advisor Jianfeng Lin for pointing out the question and helpful discussions on the proof.

\section{Criterion for nontriviality of $[t_X]_{ab}$}\label{sec2}
To prove Theorem \ref{11}, we use the following useful criterion.
\begin{Proposition}\label{21}
    Given any simply-connected smooth closed 4-manifold $X$, the following two statements are equivalent:
    \begin{enumerate}
        \item There exists a smooth $X$-bundle over a closed oriented genus $g$ surface $X\to E\to \Sigma_g$ with $w_2(TE)\neq 0$.
        \item $[t_X]_{ab}= 0\in H_1(\BDiff_\partial(X^\circ))$
    \end{enumerate}
\end{Proposition}
Note that we're working in the diffeomorphism group fixing the boundary, which justifies the following definitions.
\begin{Definition}{\ }
\begin{enumerate}
    \item Let $M$ be any smooth manifold with nonempty boundary. Let $\xi=\{M\to E\to B\}$ be a smooth fiber bundle such that $\eta=\{\partial M\to \partial E\to B\}$ is a trivial bundle. We call a trivialization $\tau:\partial E\to B\times \partial M$ of $\eta$ a boundary parameterization of $\xi$.
    \item Let $\xi_1=\{M\to E_1\to B\}$ and $\xi_2=\{M\to E_2\to B\}$ be two smooth fiber bundles with boundary parameterizations $\tau_1,\tau_2$ respectively. If there exists a bundle isomorphism $f: E_1\to E_2$ such that $\tau_2\circ f=\tau_1$, then we say $(\xi_1,\tau_1)$ is isomorphic to $(\xi_2,\tau_2)$, denoted by $(\xi_1,\tau_1)\cong(\xi_2,\tau_2)$. 
\end{enumerate}
\end{Definition}
\begin{Example}
        Let $M$ be any smooth manifold with nonempty boundary. Let $f_1,f_2\in \Diff_\partial(M)$ be two diffeomorphisms fixing $\partial M$. Consider the mapping tori $Tf_1$ and $Tf_2$ over $S^1$. By the definition of mapping torus, they can both be equipped with the canonical boundary parameterizations $\tau_c$. Then $$(Tf_1,\tau_c)\cong (Tf_2,\tau_c)$$ if and only if $$[f_1]=[f_2]\in\pi_0(\Diff_\partial(M))$$
\end{Example}
Now let's prove Proposition \ref{21}.
\begin{proof}
    \textbf{(1) implies (2):} Suppose we have a smooth fiber bundle 
    $$X\to E\to \Sigma_g$$
    with $w_2(TE)\neq 0$ Since $X$ is simply-connected, we can pick a smooth section $s$ of the projection map $\pi:E\to B=\Sigma_g$. Denote the image of $s$ by $B'$. Since $w_2(VE)\neq0$ and the normal bundle $NB'$ is isomorphic to  $VE|_{B'}$, one can use the Serre spectral sequence to show that $w_2(NB')\neq 0$. Remove a tubular neighborhood $\nu B'\subset E$ of $B'$, then we get a smooth fiber bundle $$X^\circ\to E^\circ=E-\nu B'\to B$$
    Fix any point $p\in B=\Sigma_g$ and a small open neighborhood $D\cong D^2\subset B$ around $p$. Consider the decomposition 
    $$\Sigma_g=(\Sigma_g-D)\mathop{\cup}\limits_{S^1}\bar{D}$$ 
    Since $$\Sigma_g-D\simeq \mathop{\vee}\limits_{2g} S^1$$
    the bundle structure of $(E^\circ)|_{\Sigma_g-D}$ is determined by $2g$ mapping classes $a_1,b_1,\cdots, a_g, b_g\in \pi_0(\Diff_\partial(X^\circ))$ corresponding to the standard basis of $\pi_1(\Sigma_g)$. So $$E^\circ|_{\partial\bar{D}}\cong T([a_1,b_1]\cdots[a_g,b_g])$$ is the mapping torus of $$[a_1,b_1]\cdots[a_g,b_g]\in \pi_0(\Diff_\partial(X^\circ))$$ Moreover, by the construction of mapping torus, since $[a_1,b_1]\cdots[a_g,b_g]$ fixes $\partial X^\circ$, $E^\circ|_{\partial\bar{D}}$ has a canonical trivial boundary parameterization, that is, a bundle isomorphism 
    $$\tau_c:\partial E^\circ|_{\partial\bar{D}}\cong S^1\times S^3$$
    Pick two trivializations:
    $$\tau_1:NB'|_{\Sigma^\circ}\cong \Sigma^\circ\times\mathbb{R}^4$$
    $$\tau_2:NB'|_{D^2}\cong D^2\times\mathbb{R}^4$$   
    Since $w_2(NB')\neq 0$, the clutching function $S^1\to\mathrm{SO}(4)$ must represent the generator of $\pi_1(\mathrm{SO}(4))$. $\tau_1, \tau_2$ also yield boundary parameterizations for $E^\circ|_{\Sigma^\circ}, E^\circ|_{D^2}$, still denoted by $\tau_1, \tau_2$. Since $D^2$ is contractible, $(E^\circ|_{D^2}, \tau_2)$ is isomorphic to $(D^2\times X^\circ, \tau_0)$ where $\tau_0$ is the trivial boundary parameterization given by $id_{D^2\times S^3}$. In particular, $$(E^\circ|_{\partial D^2}, \tau_2)\cong(S^1\times X^\circ,\tau_0)$$ 
    Let $\tau_0'$ is the nontrivial boundary parameterization given by:
    $$\tau_0': S^1\times S^3\to S^1\times S^3, \quad (\theta, x)\mapsto (\theta, \gamma_{\theta}(x))$$
    where $\gamma_\theta$ is a non-contractible loop in $\mathrm{SO}(4)$. Since $w_2(NB')\neq 0$, $$(E^\circ|_{\partial\Sigma_g^\circ}, \tau_1)\cong (S^1\times X^\circ, \tau_0')$$
    Note that $$(S^1\times X^\circ, \tau_0')\cong(T(t_X),\tau_c)$$
    Therefore,
    \begin{equation*}
        (T([a_1,b_1]\cdots[a_g,b_g]), \tau_c)\cong (E^\circ|_{\partial\Sigma_g^\circ}, \tau_1)\cong (S^1\times X^\circ, \tau_0')\cong (T(t_X),\tau_c)
    \end{equation*}
    Hence,
    $$[a_1,b_1]\cdots[a_g,b_g]=[t_x]\in\pi_0(\Diff_{\partial}(X^\circ))$$ 
    \textbf{(2) implies (1):} For the other direction, suppose $[t_X]_{ab}=0$, then there exists
    $$a_1,b_1,\cdots, a_g, b_g\in\pi_0(\Diff_\partial(X^\circ))$$ such that 
    $$[t_X]=[a_1,b_1]\cdots[a_g,b_g]$$
    We can use $a_1,b_1,\cdots, a_g, b_g$ to construct a smooth fiber bundle $$X^\circ\to E_1^\circ\to\Sigma_g^\circ$$ Then the restriction $$X^\circ \to E_1^\circ|_{\partial\Sigma_g^\circ}\to \partial\Sigma_g^\circ$$ is the mapping torus of $[a_1,b_1]\cdots[a_g,b_g]=[t_X]$ with canonical boundary parameterization $\tau_c$. Using the isomorphism $$(T(t_X),\tau_c)\cong (S^1\times X^\circ,\tau_0'),$$
    we can glue $(E_1^\circ,\tau_c)$ with $(E_2^\circ=D^2\times X^\circ, \tau_0)$ to get a new bundle $$X^\circ\to E^\circ\to \Sigma_g$$
    So the sphere bundles 
    $$S^3=\partial X^\circ\to E_1'\to \Sigma_g^\circ$$ and $$S^3=\partial X^\circ\to E_2'\to D^2$$ are trivialized using the boundary parameterizations $\tau_c,\tau_0$ respectively, and their clutching function turns out to be the generator of $\mathrm{SO}(4)$. Thus we may glue the $D^4$-bundle associated to the unique nontrivial rank $4$ real vector bundle $\xi$ over $\Sigma_g$ back to $X^\circ\to E^\circ\to \Sigma_g$ fiberwise along the boundary. Then we get a bundle $$X\to E\to \Sigma_g$$ with a section $s$ given by the center of $D^4$. By construction the normal bundle of the image of $s$ is isomorphic to $\xi$, which has nonzero $w_2$, so $w_2(VE)\neq 0$.
    
\end{proof}
\begin{remark}{\ }
    \begin{enumerate}
        \item This is the generalization of a criterion for triviality of $[t_X]$ in $\pi_0(\Diff_\partial(X^\circ))$ (see Proposition 2.1 in \cite{kronheimer2021dehn}). Here one just replace $\pi_1$ by $H_1$, and change the base from $S^2$ to any orientable closed surface $\Sigma_g$.
        \item Note that $TE\cong VE\oplus \pi^*T\Sigma_g$, where $VE=\mathrm{ker}(\pi_*:TE\to T\Sigma_g)$ is the vertical tangent bundle. So $w_2(TE)=0$ is equivalent to $w_2(VE)=0$, that is, the existence of a family spin structure on $VE$.
    \end{enumerate}
\end{remark}

The following two theorems are useful in the proof of Theorem \ref{11}.
\begin{Theorem}[Baraglia-Konno \cite{baraglia2022bauer}]\label{22}
    Let $(X, \mathfrak{s}_X)$ be a compact smooth $\spinc$ 4-manifold with $b_1(X) = 0$
and $b_+(X) = 3 \mod 4$. If $SW(X, \mathfrak{s}_X)$ is odd, then for any $\spinc$ family $E\to B$ with fiber $(X,\mathfrak{s}_X)$, we have $c_1(D_E)\equiv w_2(H^+)\mod 2$
\end{Theorem}

\begin{Theorem}[Baraglia-Konno \cite{baraglia2023note}]\label{23}
    For $X=K3$, the natural homomorphism $\varphi_X: \pi_0(\Diff(X))\rightarrow \Aut(H^2(X;\mathbb{Z}))$  admits a section $s_X$ over its image $\Gamma_X$ $\varphi_X$.
\end{Theorem}
\begin{remark}
    The section $s_X$ is constructed using the global Torelli Theorem of $K3$ surfaces. But for the large family of elliptic surfaces and complete intersections in \cite{baraglia2024irreducible} with nontrivial boundary Dehn twists, to the author's knowledge it's unknown whether there is a section for $\varphi_X$. So similar proof of Theorem \ref{11} doesn't work for those examples.
\end{remark}

\section{Proof of the main theorem}\label{sec3}
In this section, let $X=K3$. We are going to construct a smooth fiber bundle $X\to E\to T^2$ with $w_2(VE)\neq 0$. By Proposition \ref{21}, this implies Theorem \ref{11}. \par
By Freedman's theorem, $X$ is homeomorphic to $2(-E_8)\# 3(S^2\times S^2)$. So the Poincar\'e duals $E_i$ of diagonal classes $\Delta_i\in H_2((S^2\times S^2)_i;\mathbb{Z}),\ i=1,2,3$ form a basis of $H^2_+(X;\mathbb{Z})$.  Consider $\varphi_1,\varphi_2\in\Aut(H^2(X;\mathbb{Z}))$ given by:
$$\varphi_1(E_1)=-E_1, \varphi_1(E_2)=-E_2, \varphi_1(E_3)=E_3$$
$$\varphi_2(E_1)=E_1, \varphi_2(E_2)=-E_2, \varphi_2(E_3)=-E_3$$
and they act on the $-E_8$ components by identity.
  Let $\Gamma_X$ be the image of the natural homomorphism $\varphi_{X}: \pi_0(\Diff(X))\rightarrow \Aut(H^2(X;\mathbb{Z}))$. We know that: $$\Gamma_X=\{\phi\in\Aut(H^2(X;\mathbb{Z}))\mid \phi \text{ preserves the orientation of } H^2_+(X)\}$$ from \cite{matumoto1986diffeomorphisms}\cite{donaldson1990polynomial}\cite{friedman2013smooth}. So $\varphi_1, \varphi_2\in\Gamma_X$. By Theorem \ref{23}, we have a section $s_X:\Gamma_X\to\pi_0(\Diff(X))$ of $\varphi_X$. Let $g_1, g_2\in \Diff(X)$ with isotopy classes $ [g_1]=s(\varphi_1),[g_2]=s(\varphi_2)$, then the commutator $[g_1,g_2]$ is smoothly isotopic to identity, thus the mapping torus of $g_1, g_2$ is a smooth fiber bundle $E$ over $T^2$ with fiber $X$. \par
Let $x,y\in H^1(X;\mathbb{Z}/2)$ be a basis with $x\smile y=1$. Since $H^+(E)$ splits into the direct sum of three line bundles $L_1, L_2, L_3$ spanned by $E_1, E_2, E_3$ respectively, $w(H^+)=(1+x)(1+y)(1+x+y)=1+xy$. So $w_2(H^+)\neq 0$.
Suppose $w_2(VE)=0$, then $X\to E\to T^2$ admits a families spin structure, hence the families index $D_E\in \mathrm{KSp}(T^2)$. Therefore, $c_1(D_E)=0$, which contradicts Theorem \ref{22}. So $w_2(VE)\neq 0$ and we have proved the main theorem.
\\

\bibliographystyle{plain}
\bibliography{refs}
\end{document}